\documentclass[10pt]{article}
\usepackage{amsmath,amssymb,amsthm, graphicx,setspace,comment,mathrsfs}
\usepackage{subcaption}
\usepackage{algorithm}
\usepackage{algpseudocode}
\usepackage{caption}
\usepackage{subcaption}
\usepackage{stmaryrd}
\usepackage{tikz}
\makeatletter
\def\BState{\State\hskip-\ALG@thistlm}
\makeatother
\usepackage{dsfont}
\usepackage{mathrsfs}
\usepackage{enumerate}
\usepackage[left=3cm, right=3cm, top=2cm, bottom=2.2cm]{geometry}
\numberwithin{equation}{section} 
\numberwithin{figure}{section} 
\numberwithin{table}{section} 
\theoremstyle{theorem}
\newtheorem{thm1}{Theorem}[section]
\newtheorem{def1}[thm1]{Definition}
\newtheorem{ex1}[thm1]{Example}
\newtheorem{rm1}[thm1]{Remark}
\newtheorem{lem1}[thm1]{Lemma}
\newtheorem{prop1}[thm1]{Proposition}
\newtheorem{res1}[thm1]{Result}

\newtheorem{cor1}[thm1]{\mathcal{B}orollary}

\renewcommand{\eqref}[1]{equation (\ref{eq:#1})}
\newcommand{\eps}{\varepsilon}
\newcommand{\abs}[1]{\left|#1\right|}
\newcommand{\R}{\mathbb{R}}
\newcommand{\bs}[1]{\ensuremath{\boldsymbol{#1}}}

\def\dontshow#1{{}}
\def\R{\mathbb{R}}

\def\cond{\,\vert\,}
\def\E{\mathbb{E}}
\def\Prob{\mathbb{P}}
\def\var{\hbox{var}}
\def\mtrx#1#2{
  \left(
    \begin{array}{#1}
      #2
    \end{array}
  \right)}

\usepackage{chngcntr}
\counterwithout{figure}{section}

\title{An Additive Approximation to Multiplicative Noise}

\author{R. Nicholson$^1$ and J. P. Kaipio$^2$}

\date{%
    $^1$Department of Engineering Science, University of Auckland\\%
    $^2$Department of Mathematics, University of Auckland\\[2ex]%
    \today
}

\begin{document}
\maketitle
\begin{abstract}
Multiplicative noise models are often used instead of additive noise models in cases
in which the noise variance depends on the state.
Furthermore, when Poisson distributions with relatively small counts are approximated with 
normal distributions, multiplicative noise approximations are straightforward to implement.
There are a number of limitations in existing approaches to marginalize over 
multiplicative errors, such as positivity of the multiplicative noise term.
The focus in this paper is in large dimensional (inverse) problems for which
sampling type approaches have too high computational complexity.
In this paper, we propose an alternative approach to carry out 
approximative marginalization over the multiplicative error by embedding 
the statistics in an additive error term.
The approach is essentially a Bayesian one in that the statistics of the additive error
is induced by the statistics of the other unknowns.
As an example, we consider a deconvolution problem on random fields
with different statistics of the multiplicative noise.
Furthermore, the approach allows for correlated multiplicative noise.
We show that the proposed approach provides feasible error estimates in the sense that 
the posterior models support the actual image.
\end{abstract}

\section{Introduction}
A ubiquitous problem in science and engineering is to infer the parameter of interest, say $\bs{x}$, 
given noisy indirect measurements $\bs{y}\in\mathbb{R}^m$. 
Suppose the parameter and measurements are linked by 
a parameter-to-observable map $\bs{f}:\mathbb{R}^n\times\R^p\times\R^q\rightarrow \mathbb{R}^m$ 
\begin{align}\label{eq: general}
\bs{y}=\bs{f}(\bs{x},\bs{n},\bs{\eta})
\end{align}
where $\bs{y}\in\R^m$ is the (observation) data,
$\bs{x}\in\R^n$ is the primary (interesting) unknown and
$\bs{n}\in\R^p$ and $\bs{\eta}\in\R^q$ denote uninteresting related random variables
which can often be interpreted as noise. 
The first task would then be to marginalize over the uninteresting variables.
In the context of inverse problems which are the focus in this paper, we often have $n\ge m$.

The most common model for $\bs{f}(\bs{x},\bs{n},\bs{\eta})$ is the additive error model \cite{Tarantola2004,kaipio2005}
\begin{align}\label{eq: additive}
\bs{y}=\bs{A}(\bs{x})+\bs{\eta}
\end{align}
where the mapping $\bs{A}: \bs{x}\mapsto\bs{y}$ is referred to as the forward map (problem).
However, in several imaging modalities including optical coherence tomography (OCT) \cite{Wong2010,Yin2013}, 
ultrasound \cite{Burckhardt1978,Michailovich2006}, synthetic aperture radar (SAR) imaging \cite{Foucher2001,Tison2004}, 
and electrical impedance tomography (EIT) \cite{Borcea1999,Zhang2015}, noise can be proportional to the data. 
In such a case, we have 
\begin{align}\label{eq: noadd}
\bs{y}=\bs{n}\odot\bs{A}\left(\bs{x}\right)
\end{align}
where $\odot$ denotes component-wise (Hadamard) product. 
Such a set up is usually referred to as the multiplicative noise model. 
Moreover, in many of these cases there may simultaneously be additive noise present, 
see, for example, \cite{Krissian2007,Durand2010,kaipio2005,Garcia2011}, so that we can write
\begin{align}\label{eq: fullmodel}
\bs{y}=\bs{n}\odot\bs{A}\left(\bs{x}\right)+\bs{\eta}.
\end{align}
In most papers, see for example \cite{Krissian2007,Durand2010}, the effects of any additive errors $\bs{\eta}$ have been 
assumed to be small compared to 
the effects of the multiplicative noise $\bs{n}$, and thus the additive error term 
has often been neglected. 
Furthermore, the multiplicative noise has systematically been assumed to be
mutually independent.
In the current paper however, we retain the additive error term. 

In this paper, we take a model discrepancy style approach to transform Equation (\ref{eq: fullmodel}) to 
Equation (\ref{eq: additive}) with a modified additive error term,
into which both the additive and multiplicative errors are embedded. 
The approach is based on a joint model $\pi(\bs{x},\bs{n},\bs{\eta})$ 
and the computation of the approximate statistics
of the adjusted additive error term, followed by approximate marginalization.
This procedure yields an approximate posterior model $\pi( \bs{x} \cond \bs{y})$.

The paper is organised as follows. In Section \ref{sec: Noise}, we review the marginalization of
noise terms 
in the Bayesian framework.
In Section~\ref{sec: StdApp}, we give a brief review of 
the methods used to deal with multiplicative noise.
Section~\ref{sec: BAE} outlines the approximation of the noise statistics and the subsequent marginalization,
which approach is sometimes referred to as the Bayesian approximation error (BAE) approach 
\cite{kaipio2005,Kaipio2007,Kaipio2013}.
The multiplicative noise term is not assumed to be uncorrelated.
In Section~\ref{sec:DeblurringApplication}, we consider a deconvolution example with different distributions for the 
multiplicative noise term, including correlated multiplicative noise. 
\dontshow{The results are compared to those based on the same transformation without updating the additive error term so as 
to  take into account the effects of the multiplicative noise.}

\section{Exact marginalization over additive and multiplicative terms}
\label{sec: Noise}
In this paper, we assume that  the noise terms $\bs{n}$ and $\bs{\eta}$ and the parameter of interest $\bs{x}$ 
are pair-wise mutually independent. 
Thus, joint model of the noise terms and the parameter can be stated as 
$\pi(\bs{x},\bs{n},\bs{\eta})=\pi_x(\bs{x})\pi_n(\bs{n})\pi_{\eta}(\bs{\eta})$.
Furthermore, in line with the literature \cite{Shi2008,Aubert2008,Durand2010,Steidl2010,Zhao2014}, 
we assume that the multiplicative noise is i.i.d., so that $\pi_n(\bs{n})=\prod_{i=1}^m\pi_{n_i}(n_i)$. 

The likelihood is obtained formally by marginalization
\begin{align}
\label{eq: fulllike}
\pi(\bs{y}|\bs{x})=\iint\pi(\bs{y}|\bs{x},\bs{n},\bs{\eta})\pi_n(\bs{n})\pi_{\eta}(\bs{\eta})\;d\bs{n}d\bs{\eta}=\iint \delta(\bs{y}-\bs{n}\odot\bs{A}(\bs{x})-\bs{\eta})\pi_n(\bs{n})\pi_{\eta}(\bs{\eta})\;d\bs{n}d\bs{\eta}.
\end{align}
where $\delta(\cdot)$ is the Dirac distribution.
We now look at the three individual cases of interest.

In the purely additive noise model, we set $\pi(\bs{n})=\delta(\bs{n}-\bs{1})=\prod_{i=1}^m\delta(n_i-1)$. 
Thus, (\ref{eq: fulllike}) can be written as
\begin{align}
\pi(\bs{y}|\bs{x})&=\iint \delta(\bs{y}-\bs{n}\odot\bs{A}(\bs{x})-\bs{\eta})\delta(\bs{n}-\bs{1})\;d\bs{n}\pi_{\eta}(\bs{\eta})\;d\bs{\eta}
=\int\delta(\bs{y}-\bs{A}(\bs{x})-\bs{\eta})\pi_{\eta}(\bs{\eta})\;d\bs{\eta}\nonumber\\
&=\pi_\eta(\bs{y}-\bs{A}(\bs{x})),
\end{align}

In the purely multiplicative noise model, we set $\pi(\bs{\eta})=\delta(\bs{\eta})=\prod_{i=1}^m\delta(\eta_i)$, then
(\ref{eq: fulllike}) can be written as
\begin{align}\pi(\bs{y}|\bs{x})&=\iint \delta(\bs{y}-\bs{n}\odot\bs{A}(\bs{x})-\bs{\eta})\delta(\bs{\eta})\;d\bs{\eta}\pi_{n}(\bs{n})\;d\bs{n}
=\int\delta(\bs{y}-\bs{n}\odot\bs{A}(\bs{x}))\pi_{n}(\bs{n})\;d\bs{n}\nonumber\\
&=\prod_{i=1}^m\left(\frac{1}{\abs{A_i(\bs{x})}}\pi_{n_i}\left(\frac{y_i}{A_i(\bs{x})}\right)\right),
\end{align}
where $A_i(\bs{x})$ is the $i$th component of $\bs{A}(\bs{x})$. 

In the case of simultaneous multiplicative and additive noise terms, due to Fubini's theorem, 
the integrations in (\ref{eq: fulllike})  can be carried out in either order 
resulting in either
\begin{align}\label{eq: int1}
\pi(\bs{y}|\bs{x})&=\iint \delta(\bs{y}-\bs{n}\odot\bs{A}(\bs{x})-\bs{\eta})\pi_n(\bs{n})\pi_{\eta}(\bs{\eta})\;d\bs{n}d\bs{\eta}\nonumber\\&=\int\pi_n(\bs{n})\pi_{\eta}(\bs{y}-\bs{n}\odot\bs{A}(\bs{x}))\;d\bs{n}.
\end{align}
or
\begin{align}\label{eq: int2}
\pi(\bs{y}|\bs{x})&=\iint \delta(\bs{y}-\bs{n}\odot\bs{A}(\bs{x})-\bs{\eta})\pi_n(\bs{n})\pi_{\eta}(\bs{\eta})\;d\bs{\eta}d\bs{n}\nonumber\\&=\prod_{i=1}^m\left(\frac{1}{\abs{A_i(\bs{x})}}\right)
\int\prod_{i=1}^m\left(\pi_{n_i}\left(\frac{y_i-\eta_i}{A_i(\bs{x})}\right)\right)\pi_\eta(\bs{\eta})\;d\bs{\eta}.
\end{align}
Unfortunately, the integrals defined as in either (\ref{eq: int1}) or (\ref{eq: int2}) 
cannot be computed analytically for general multiplicative noise models $\pi_{\eta}(\bs{\eta})$. 

\section{Approaches to handle multiplicative noise models}
\label{sec: StdApp}
For the remainder of the paper we will consider linear forward models 
$\bs{A}(\bs{x})=\bs{A}\bs{x}$, as is the case in deblurring (setting $\bs{A}=\bs{I}$ is the case of denoising). 
Furthermore, since the focus of the present paper is in inverse problems and since some approaches depend
directly on properties of the unknown (such as positivity), we refer directly on posterior models.
Moreover, since the proposed approach is targeted on relatively large dimensional problems, we
will consider the computation of MAP estimates and the Laplace approximations for 
the posterior covariances only.

There are several approaches documented in the literature for dealing with multiplicative noise.
 Many of the techniques are framed in the context of denoising. 
 Moreover, it is often assumed that $\bs{x}\ge\bs{0}$ and has bounded variations 
 (i.e. $\bs{x}\in BV(\Omega)$) and thus the total variation (TV) prior is used, see for example
  \cite{Shi2008,Aubert2008,Durand2010,Steidl2010,Rodriguez2013,Zhao2014}.  In the approach proposed below, however, we do not need to assume positivity or boundedness
of the primary unknown $x$.

{\bf The} \bs{\log} {\bf model (multiplicative noise only):} 
The most common of these techniques is to simply apply 
the logarithm transform, resulting in a problem of the form of (\ref{eq: additive}), see for example \cite{Guo1994,Durand2010}. 
However, there are some drawbacks to applying the logarithm transform method. 
Firstly, if any of the components of the data $\bs{y}$, the model prediction $\bs{A}\bs{x}$ 
or the noise term $\bs{n}$ are negative, the method fails. Secondly if one in fact retains 
the additive error, as in Equation (\ref{eq: fullmodel}), the logarithm transform is of little use. 
Thirdly, it has been noted that one cannot directly apply standard additive noise removal algorithms 
and that such a method does not produce satisfactory results \cite{Aubert2008}. 
Such an approach leads 
to the following MAP estimate for a general prior on $\pi_x(\bs{x})$ 
\dontshow{placed on $\bs{x}$,}
\begin{align}\label{eq: logMAP}
\bs{x}_{\rm MAP}=\max_{\bs{x}}\pi_\xi\left(\log(\bs{y})-\log\left(\bs{A}\bs{x}\right)\right)\pi_x(\bs{x}),
\end{align}
where $\pi_\xi$ is the density of $\xi=\log(\bs{n})$. 
Iterative method can then be used to solve for $\bs{x}_{\rm MAP}$.
The basic idea of transforming multiplicative noise to additive noise is in principle
similar  to the procedure 
we propose in the current paper, except that, in this paper, the measurements are not transformed. 

 {\bf The AA model:} The so-called AA model was derived in \cite{Aubert2008} for the MAP estimate 
 under the assumption that the multiplicative noise follows a Gamma distribution and under the prior assumption 
 that $\bs{x}$ has bounded variations and is positive \cite{Aubert2008}. 
 The MAP estimate is then 
 computed for a forward map
 $\bs{A}$
 \begin{align}\label{eq: AAMAP}
\bs{x}_{\rm MAP}=\min_{\bs{x}}\sum_{i=1}^n\left(L\left(\log(A_i\bs{x})+\frac{y_i}{A_i\bs{x}}\right)+\gamma\phi(x_i)\right),
\end{align}
where $A_i$ denotes the $i$th row of $\bs{A}$ and the term $\gamma\phi(x_i)$ is 
induced by the total variation prior on $\bs{x}$. 
The computation of the MAP estimate in this case also requires iterative methods even when 
the prior on $\bs{x}$ were Gaussian.
Furthermore, the likelihood potential
is not always strictly convex although the existence of a minimiser was proven in
\cite{Aubert2008}.

 {\bf The separable model:}  The separable model was introduced in \cite{huang2013} and takes into account both additive and multiplicative noise. Furthermore, for several different multiplicative noise models, closed form functionals are derived to find the respective MAP estimates. Here we give a brief outline of how  one can derive the posterior. In accordance with \cite{huang2013}, The separable model takes the form 
\begin{align}\label{eq: sep}
\bs{y}=\bs{f}(\bs{x},\bs{n},\bs{\eta})=\bs{n}\odot\left(\bs{A}\left(\bs{x}\right)+\bs{\eta}\right).
\end{align}
The key ingredient to dealing with the separable model is the introduction of an intermediate variable, 
\begin{align}
\bs{u}=\bs{A}\bs{x}+\bs{\eta}.
\end{align}
Introduction of this intermediate variable results in the posterior of interest being given by
  \begin{align}\label{eq: CRposts}
\hat{\pi}_{\rm post}(\bs{u},\bs{x}|\bs{y})\propto\pi(\bs{u},\bs{x})\pi(\bs{y}|\bs{u},\bs{x})=\pi_x(\bs{x})\pi(\bs{u}|\bs{x})\pi(\bs{y}|\bs{u}).
\end{align}
Each of the conditional densities in (\ref{eq: CRposts}) can be derived similarly to how the likelihood densities were dealt with in Section \ref{sec: Noise}. Firstly,
 \begin{align}
\pi(\bs{y}|\bs{u})&=\int \delta(\bs{y}-\bs{n}\odot\bs{u})\pi_n(\bs{n})\;d\bs{n}\nonumber\\&=\prod_{i=1}^m\left(\frac{1}{\abs{u_i}}\pi_{n_i}\left(\frac{y_i}{u_i}\right)\right),
\end{align}
and secondly,
 \begin{align}
\pi(\bs{u}|\bs{x})&=\int \delta(\bs{\eta}-\bs{u}-\bs{A}\bs{x})\pi_\eta(\bs{\eta})\;d\bs{\eta}\nonumber\\&=\pi_\eta\left(\bs{u}-\bs{A}\bs{x}\right),\end{align}
hence the posterior can be written as
 \begin{align}
\pi(\bs{u},\bs{x}|\bs{y})=\prod_{i=1}^m\left(\frac{1}{\abs{u_i}}\pi_{n_i}\left(\frac{y_i}{u_i}\right)\right)\pi_\eta\left(\bs{u}-\bs{A}\bs{x}\right)\pi_x(\bs{x}).
\end{align}
The MAP estimate, $(\bs{u}_{\rm MAP},\bs{x}_{\rm MAP})$ is shown in closed form for several prior densities on the multiplicative noise in \cite{huang2013}. The main drawback to this method is the lack of convexity for the functionals which need to be minimised in order to calculate the MAP when $\bs{u}$ is not strictly positive. Furthermore the computation of the MAP estimate, 
again requiring iterative method irrespective of the prior model, 
will need to be carried out for the number of primary unknowns and the number of measurements. 

Other methods for dealing with multiplicative noise in the denoising context include filtering type approaches such as those discussed in \cite{Garcia2011} and the use and construction of similarity measures \cite{Teuber2012}. For filtering type methods the problem is framed in the so called state-space formalism. On the other hand, approaches using similarity measures usually set values of the restored image to some weighted mean of the surrounding pixels, where the weights depend on the similarity of the pixels.

\section{Approximate marginalization of multiplicative noise}
\label{sec: BAE}
In this paper, we carry out approximative marginalization over both the additive and multiplicative noise terms.
In the inverse problems literature, this approach is referred to as the 
Bayesian approximation error  approach (BAE) since the approximative marginalization is carried over the prior distribution.
The BAE was introduced in \cite{kaipio2005,Kaipio2007} to take into account 
the discrepancy between accurate and reduced order models. 
Since that, the approach has been extended, for example, 
to account for errors and uncertainties related to
uninteresting distributed parameters in PDE's \cite{kolehmainen2011},
errors in the geometry of the domain \cite{nissinen2011a},
unknown boundary data \cite{lehikoinen2007},
approximation of the (physical) forward map \cite{tarvainen2010},
and state estimation problems \cite{jhuttunen2007a,lehikoinen2009,Lipponen2010}.
For a more general discussion of the approach and a more extended 
list of extensions, see \cite{Kaipio2013}.
Below, we adapt the approach to the context of multiplicative noise. 

The goal is to embed the  additive and multiplicative noise terms
in an {\it approximate} additive error only model.
With such an approximation together with linear forward and normal prior models,
the computation of the approximate MAP estimate and the approximate 
posterior covariance reduces to linear algebra.

With the present observation model, we can write
\begin{align}
\bs{y}&=\bs{n}\odot\bs{A}\bs{x}+\bs{\eta}\nonumber \\
&=\bs{A}\bs{x}+\left(\bs{n}-\bs{1}\right)\odot\bs{A}\bs{x}+\bs{\eta}\nonumber \\
&=\bs{A}\bs{x}+\bs{\eps}+\bs{\eta}\nonumber\\
&=\bs{A}\bs{x}+\bs{e},
\end{align} 
which is alternative and {\it exact} additive error model to use in place of (\ref{eq: fullmodel}).  
Exact marginalization over $\bs{e}$ would then yield
the likelihood model $\pi(\bs{y}\cond \bs{x}) = \pi_{e\cond x}(\bs{y} -\bs{A}\bs{x}\cond\bs{x})$
the computation of which is, however, not generally possible analytically.
In the BAE approach, at this stage,
one (usually) makes the normal approximation 
$\pi_{e\vert x}(\bs{e})=\mathcal{N}(\bs{e}_{*\vert x},\bs{\Gamma}_{e\vert x})$, see, 
for example, \cite{kaipio2005,Kaipio2007,Kaipio2013}.
We note, however, that some work has been carried out on retaining the full density of the errors, 
$\pi_e(\bs{e})$ \cite{Calvetti2014,Calvetti2017}.
Furthermore, in theory, the full density of the approximation errors could be calculated as 
the product density of $\bs{p}\odot\bs{q}$ for $\bs{p}=\bs{n}-\bs{1}$ and $\bs{q}=\bs{A}\bs{x}$, 
see for example \cite{Rohatgi2015}. 

\dontshow{
In the BAE approach, a further (technical) approximation 
$\pi_{e\cond x}(y - Ax\cond x)\approx \pi_{e}(y - Ax)$ is sometimes carried out, see for example
\cite{kaipio2005}.
This approximation leads to the so-called enhanced error model and it is often 
computationally less heavy to determine than the full BAE model.
For the obvious reason that $e = e(x)$, this further approximation is not always a feasible one,
see for example \cite{Kaipio2007}.
Below, however, we employ the enhanced error model, that is, we approximate
$\bs{e}_{*\vert x} \approx \bs{e}_{*}$ and $\bs{\Gamma}_{e\vert x} \approx \bs{\Gamma}_{e}$.
Below, we make the standard assumption that $(n,\eta,x)$ are mutually independent.
}

For the mean $\E(\bs{e}) = \bs{e}_\ast$, we get
\begin{align}
\bs{e}_*=\bs{\eta}_*+(\bs{n}_*-\bs{1})\odot\bs{A}\bs{x}_*.
\end{align}
In case we have, as is the standard assumption, $\E(\bs{n}) = \bs{n}_* = \bs{1}$ and   $\E(\bs{\eta}) = \bs{0}$, 
we also have $\bs{e}_*=\bs{0}$.
For the joint covariance matrix 
\[
\bs{\Gamma}_{x,e} = \mtrx{cc}{\bs{\Gamma}_{xx} & \bs{\Gamma}_{xe} \\ \bs{\Gamma}_{ex} & \bs{\Gamma}_{ee} }
\]  we have
\begin{eqnarray}
\label{eq:ApprCov}
\bs{\Gamma}_{ee} &=& \bs{\Gamma}_{\eta\eta}+\bs{\Gamma}_{nn}\odot\bs{A}\bs{\Gamma}_{xx}\bs{A}^T \\
\bs{\Gamma}_{ex} &=& \bs{\Gamma}_{\eta x}
\end{eqnarray}
due to the assumption that $n$ is uncorrelated with both $(x,\eta)$.
Note that we do not have to assume the uncorrelatedness of $n$ or mutual
uncorrelatedness  of $(x,\eta)$.
For the conditional covariance $\Gamma_{e\vert x}$, we have
\[
\bs{\Gamma}_{e\vert x}
=  \bs{\Gamma}_{\eta\eta} + \bs{\Gamma}_{nn}\odot\bs{A}\bs{\Gamma}_{xx}\bs{A}^T 
    - \bs{\Gamma}_{\eta x}\bs{\Gamma}_{xx}^{-1}\bs{\Gamma}_{x\eta}
\]
If $\bs{\Gamma}_{e\vert x}$ has full rank, the approximate likelihood model can then
be written as
\[
\pi(\bs{y}\cond \bs{x}) = \pi_{e\vert x}(\bs{y} - \bs{A}\bs{x}\cond \bs{x}) \propto \exp\left( -\frac12  \left\Vert 
\bs{L}_{e\vert x}\left( \bs{y} - \bs{A}\bs{x} - \bs{\eta}_\ast - \bs{\Gamma}_{\eta x}\bs{\Gamma}_{xx}^{-1}(\bs{x} - \bs{x}_\ast) \right)   
\right\Vert_2^2 \right)
\]
where $\bs{L}_{e\vert x}^T \bs{L}_{e\vert x} = \bs{\Gamma}_{e\vert x}^{-1}$.
However, with the typical assumptions of  i.i.d. additive and multiplicative noise and 
mutually uncorrelatedness of $x$ and $\eta$, 
we have
\begin{eqnarray}
\bs{\Gamma}_{ee} %&=& \sigma^2_\eta\bs{I}+\sigma_n^2\bs{I}\odot\bs{A}\bs{\Gamma}_x\bs{A}^T \\
%  &=& \sigma^2_\eta\bs{I}+\sigma_n^2\bs{I}\odot\text{diag}\left(\bs{A}\bs{\Gamma}_x\bs{A}^T\right)
  &=& \sigma^2_\eta\bs{I}+\sigma_n^2\text{diag}\left(\bs{A}\bs{\Gamma}_x\bs{A}^T\right)
\end{eqnarray}
so that we have 
$\bs{\Gamma}_{e\vert x} = \bs{\Gamma}_{ee}$.
Again, if $\bs{\Gamma}_{ee}$ is full rank, we can write
$\bs{L}_e^T \bs{L}_e = \bs{\Gamma}_{ee}^{-1}$ which results in the approximate likelihood model
\[
\pi(\bs{y}\cond \bs{x}) \propto \exp\left(-\frac12 \Vert \bs{L}_e(\bs{y} - \bs{A}\bs{x}) \Vert_2^2 \right)
\]
Note that the structure of the covariance $\Gamma_{ee}$ is, in the general case, nontrivial
and depends also on the prior covariance $\Gamma_{xx}$, which is the case in BAE type approaches.
As far as the authors are aware, correlated multiplicative noise has not previously been considered 
in the literature.
%%%%%%%%%%

\section{Application to deblurring}
\label{sec:DeblurringApplication}
We consider an image deblurring (deconvolution) example with three different multiplicative noise statistics:
normal, Gamma and uniform distributions.
The image used assumes both positive and negative values.
For each case additive noise of level with standard deviation corresponding to
1\% of the range of the noiseless observations.
We also consider the same example with normal multiplicative noise that is spatially correlated.
Since the focus in this paper is on the multiplicative noise, we take the image to be uncorrelated with the 
additive noise component.

%%%%%%%%%%%%%%%%%%%%%%%
\subsection{Multiplicative and additive noise models}
Without loss of generality, we set $\mathbb{E}(\bs{n})=\bs{1}$ for all cases 
as is customary \cite{Aubert2008,Shi2008,Li2010}. 
In this section, we take the components of $n$ to be iid so that also $\Gamma_{nn} = \sigma_n^2 I$.
We also take the additive noise model to be $\bs{\eta}\sim\mathcal{N}(\bs{0},\sigma^2_\eps\bs{I})$.
A correlated additive noise model is straightforward to handle as in Section~\ref{sec: BAE}.
We consider three different distributions for the multiplicative noise $\bs{n}$ and scale
them so that the variances coincide.
Furthermore, in two cases, the probability $\Prob(n_i)<0$ does not vanish.

The first model for $\bs{n}$ is the iid Gamma distribution which has been 
the most common model for multiplicative noise
 \cite{Aubert2008,Shi2008,Durand2010},
\begin{align}
n_i\sim\Gamma(\alpha,\beta),\quad i=1,2,\dots, m
\end{align}
where $\alpha$ and $\beta$ are the shape and scale parameters, respectively. 
For $\bs{n}$ such that is $\E(\bs{n}) = \bs{1}$ we can also write
\begin{align}
n_i\sim\Gamma\left(L,\frac{1}{L}\right),\quad i=1,2,\dots, m.
\end{align}
We set $L=1$ so that $\var(n_i) = 1$.
For the Gamma distribution, $\Prob(n_i<0) = 0$.

The second model for $\bs{n}$ which is less seldom considered is the iid normal model
\begin{align}
n_i\sim\mathcal{N}(1,\sigma^2),\quad i=1,2,\dots, m,
\end{align}
where throughout the literature $\sigma\leq 0.2$ is referred to as {\it tiny noise},
see, for example, \cite{Aubert2008}. 
The assumption of tiny noise is done in an attempt to avoid multiplicative noise terms becoming negative, 
as discussed in more detail in Section~\ref{sec: StdApp}. 
In the approach proposed in this paper, we do not need to make such an assumption and we set, 
again, $\sigma = 1$
which results in $\Prob(n_i<0)\approx 0.15$.

As the third model, we consider multiplicative noise with iid uniform distribution,
\begin{align}
n_i\sim\mathcal{U}(1-\nu,1+\nu),\quad i=1,2,\dots, m
\end{align}
and we set $\nu = \sqrt{3}$ so that, again, $\var(n_1) = 1$ and
which results in $\Prob(n_i<0)\approx 0.21$.

\begin{figure}[h!]
\includegraphics[width=16cm]{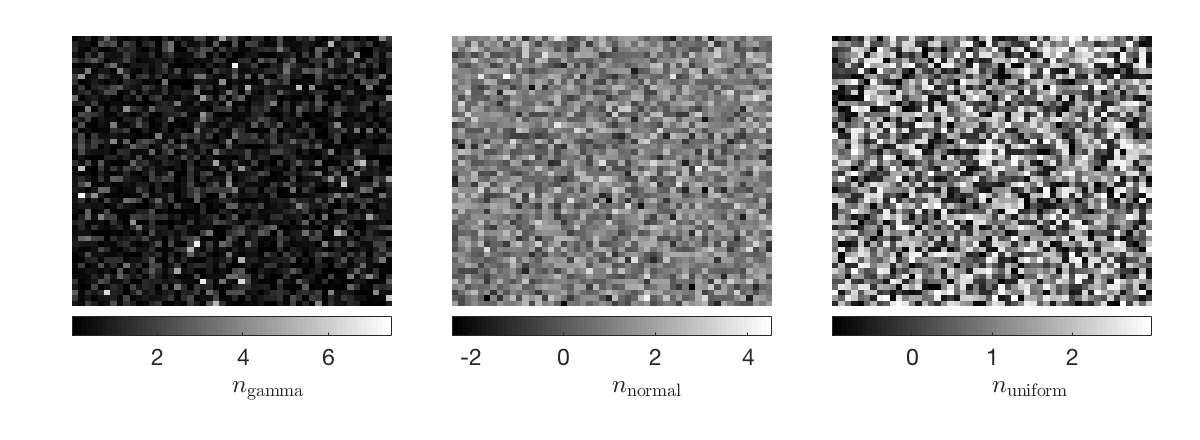} 
\caption{Draws from the different iid multiplicative noise distributions. 
Left: Gamma distribution $\Gamma(1,1)$. 
Centre:  normal distribution $\mathcal{N}(\bs{1},\bs{I})$. 
Right: uniform distribution $\mathcal{U}(1-\sqrt{3},1+\sqrt{3})$.}
\label{fig: noisedraws}
\end{figure}

Draws from the three multiplicative noise models are shown in Fig.~\ref{fig: noisedraws}. 
%%%%

%%%%
\begin{figure}[h!]
\includegraphics[width=16cm]{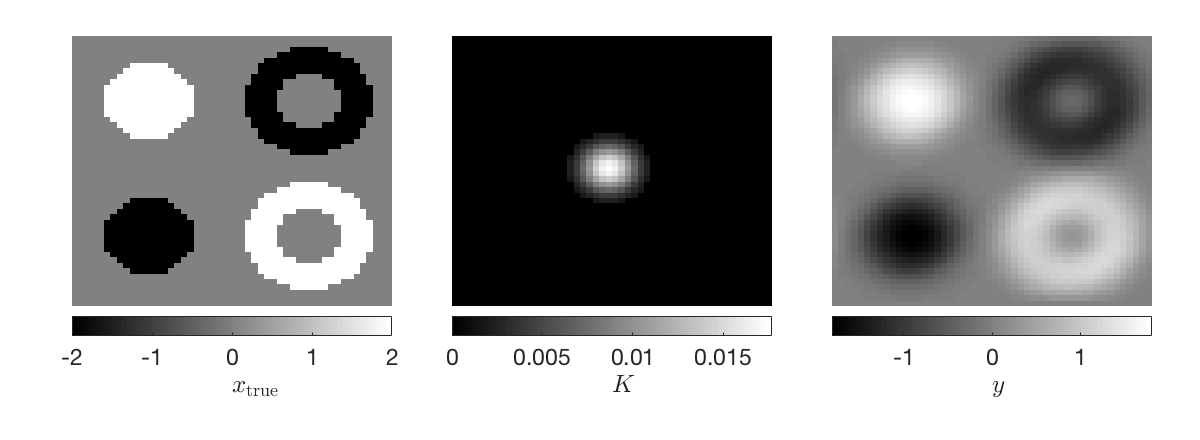} 
\caption{Left: the target image $\bs{x}_{\rm true}$. Centre: the Gaussian convolution kernel $\mathcal{K}$
centred at the center of the image. 
Right: the  blurred noiseless image $\bs{K}*\bs{x}$.}\label{fig: presetup}
\end{figure}
%%%%%%%%%%%%%%%%%%%

\subsection{The target, the observations and the prior model}
\begin{figure}[h!]
\includegraphics[width=16cm]{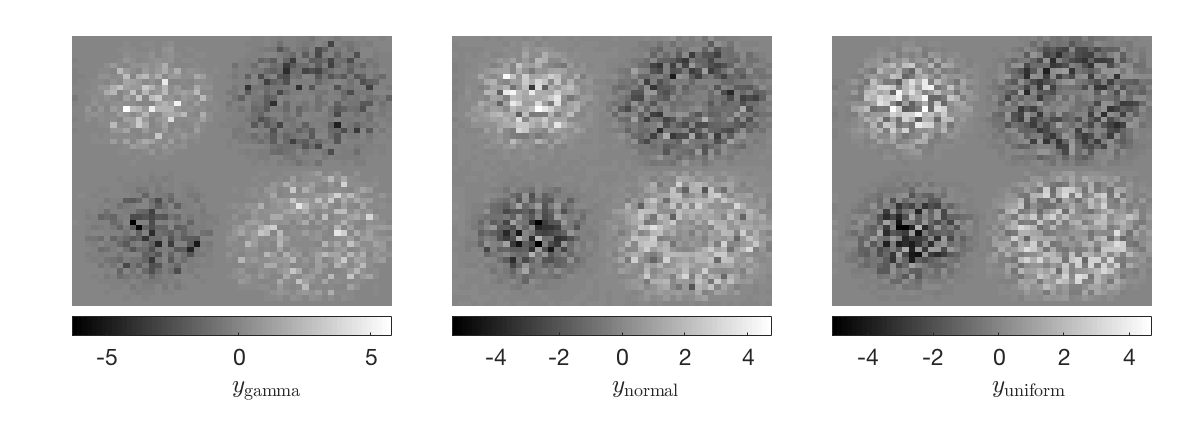} 
\caption{Data corrupted by multiplicative noise generated from left: Gamma, centre: Normal, and right: Uniform distributions.}\label{fig: data}
\end{figure}

For all examples, we specify a $50\times50$ pixel target image shown in Fig.~\ref{fig: presetup}. 
We blur the image with a symmetric Gaussian blurring kernel
\begin{align}
\mathcal{K}(s_1,s_2) = \frac{1}{2\pi\kappa^2}\exp\left(-\frac{s_1^2+s_2^2}{2\kappa^2}\right),
\end{align}
with $\kappa=5$ also shown in Fig.~\ref{fig: presetup}.
Both the image and the kernel are taken to be piecewise constant in a grid
with rectangular elements.
We take the forward operator to be the circulant convolution operator $\mathcal{K}$ \cite{Calvetti2005}
\begin{align}
\bs{y}=\bs{n}\odot \left(\bs{K}*\bs{x}\right)+\bs{\eta}=\bs{n}\odot \bs{A}\bs{x}+\bs{\eta},
\end{align}
where  $\bs{K}$ is the circulant realization of the kernel $\mathcal{K}$ 
and, further,  $\bs{A}$ is the realization of $\bs{K}$ in matrix form. 
The blurred (noiseless) image $\bs{y}=\bs{A}\bs{x}$ is also shown in Fig.~\ref{fig: presetup}.
The observations with the three different multiplicative noise models are shown in Fig.~\ref{fig: data}.

\begin{figure}[t!]
\includegraphics[width=16cm]{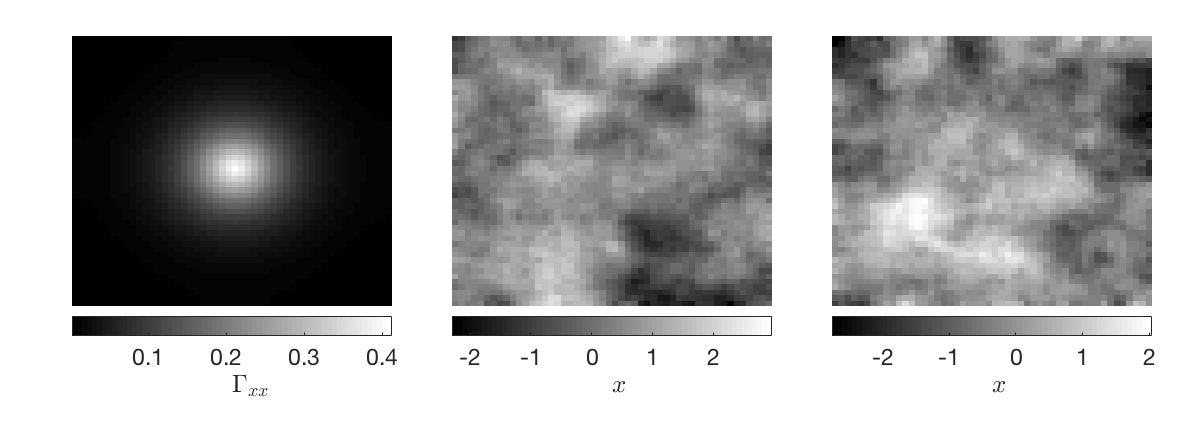} 
\caption{Left: The correlation function induced by the PDE based model with 
$c_1=10^{-1}$ and $c_2=20$.
Center and right: two draws from the prior model.}\label{fig:PriorCov}
\end{figure}

%%%%%%%%%%%%%%%% Prior model
In this paper, we employ a normal prior model $\bs{x}\sim\mathcal{N}(\bs{x}_*,\bs{\Gamma}_x)$. 
The mean  of $\bs{x}$ is set to be spatially homogeneous $\E(\bs{x}) = \bs{x}_* = x_*\bs{1}$. 
For the prior covariance matrix, we employ  so-called PDE-based covariance matrices 
\cite{Stuart2010,Bui-Thanh2013}. 
More specifically, we take
\begin{align}\label{eq: pCov}
\bs{\Gamma}_{xx}=\left(c_1\left(c_2\bs{G}+\bs{M}\right)\right)^{-2}  = \left(L_x^T L_x  \right)^{-1},
\end{align}
where $c_1$ is a constant inversely proportional to the variance, $c_2$ is a constant which 
controls the correlation length.
The matrix square root $L_x$ is the whitening operator of the random field
\cite{RueBook}.
The matrices $\bs{G}$ and $\bs{M}$  are the stiffness and mass matrices, respectively
\begin{align}
G_{ij}=\int_\Omega\nabla\phi_i\cdot\nabla\phi_j\;d\bs{s} \quad M_{ij}=\int_\Omega\phi_i\phi_j\;d\bs{s},\quad i,j=1,2,\dots,n.
\end{align}
The parameters are set as $c_1=10^{-1}$ and $c_2=20$ so that the range of $x$ and the correlation length
are approximatively consistent with the structure of the target image, and we also set  $\bs{x}_*=\bs{0}$, see Fig.~\ref{fig:PriorCov}
for the covariance function and two draws from the prior model.

%%%%%%%%%%%%%%%%%

\subsection{Reconstructions with spatially uncorrelated multiplicative noise}
\begin{figure}[h!]
\includegraphics[width=16cm]{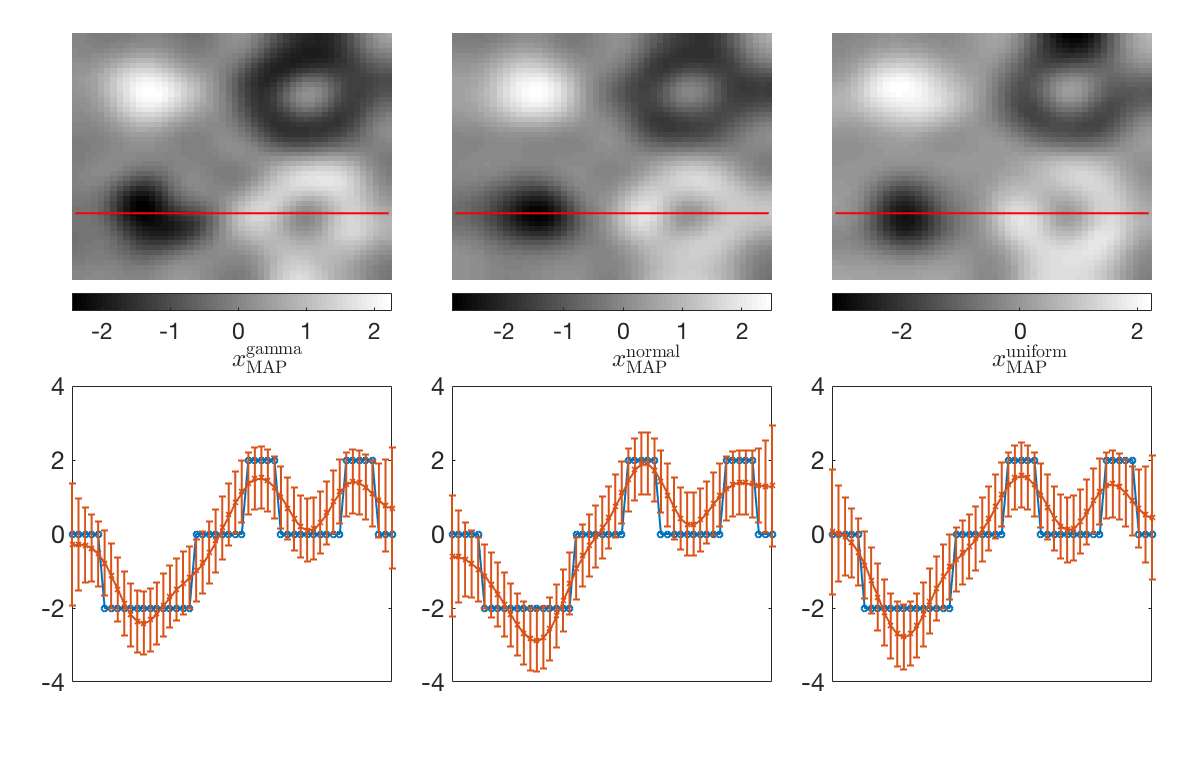} 
\caption{Top row.
The MAP estimates attained by using the BAE approach with different iid multiplicative noise models.
Left: Gamma, Centre: normal and Right: uniform noise models. 
Bottom row.
Cross sections of the actual target and reconstructions with  approximate MAP $\pm 3\sqrt{\bs{\Gamma}_{x\vert y}(k,k)}$
intervals along the lines in the top row reconstructions.}\label{fig: baeMAPS}
\end{figure}

The reconstructions computed using the proposed approximation, 
denoted $\bs{x}_{\rm MAP}^{\rm gamma}$, $\bs{x}_{\rm MAP}^{\rm normal}$, 
and $\bs{x}_{\rm MAP}^{\rm uniform}$
are shown in Fig.~\ref{fig: baeMAPS}. 
Furthermore, in the bottom row of Figure \ref{fig: baeMAPS} we show the estimates and 
posterior confidence intervals along the cross section shown in the images in the top row. 
We see that embedding the multiplicative noise into the additive error leads to feasible results
in the sense that the actual target is supported by the approximative MAP $ \pm 3\sqrt{\bs{\Gamma}_{x\vert y}(k,k)}$
intervals. 
It is noteworthy that the fact that, in the case of normal and uniform multiplicative noise distributions,
$n$ exhibits negative samples.
Clearly, this does not constitute a problem for the proposed approach.
The feasibility of the posterior estimates is similar with all three distributions of the multiplicative noise.
The fact that the estimates obtained here are fairly smooth in comparison to the true image is due to the use 
of a Gaussian smoothness prior. 
To reconstruct the sharp edges one could employ a TV type prior.
Even with a normal approximation for the posterior, 
this would result in the need for iterative methods to compute the MAP estimate.

%%%%

\subsection{Reconstructions with spatially correlated multiplicative noise}
The derivations of the existing methods to handle multiplicative noise are largely based on the assumed iid property.
In the proposed approach, such an assumption does not need to be done, as indicated by the approximate 
joint covariance $\bs{\Gamma}_{e,x}$ in 
Section~\ref{sec: BAE}.
With a deblurring problems such as the present example, it is clear that when the spatial correlation
structure of the multiplicative noise gets more complicated, we can expect the actual estimation errors to increase.
This can be expected, in particular, with noise distributions with positive spatial and
increasing correlation length.

\begin{figure}[h!]
\includegraphics[width=16cm]{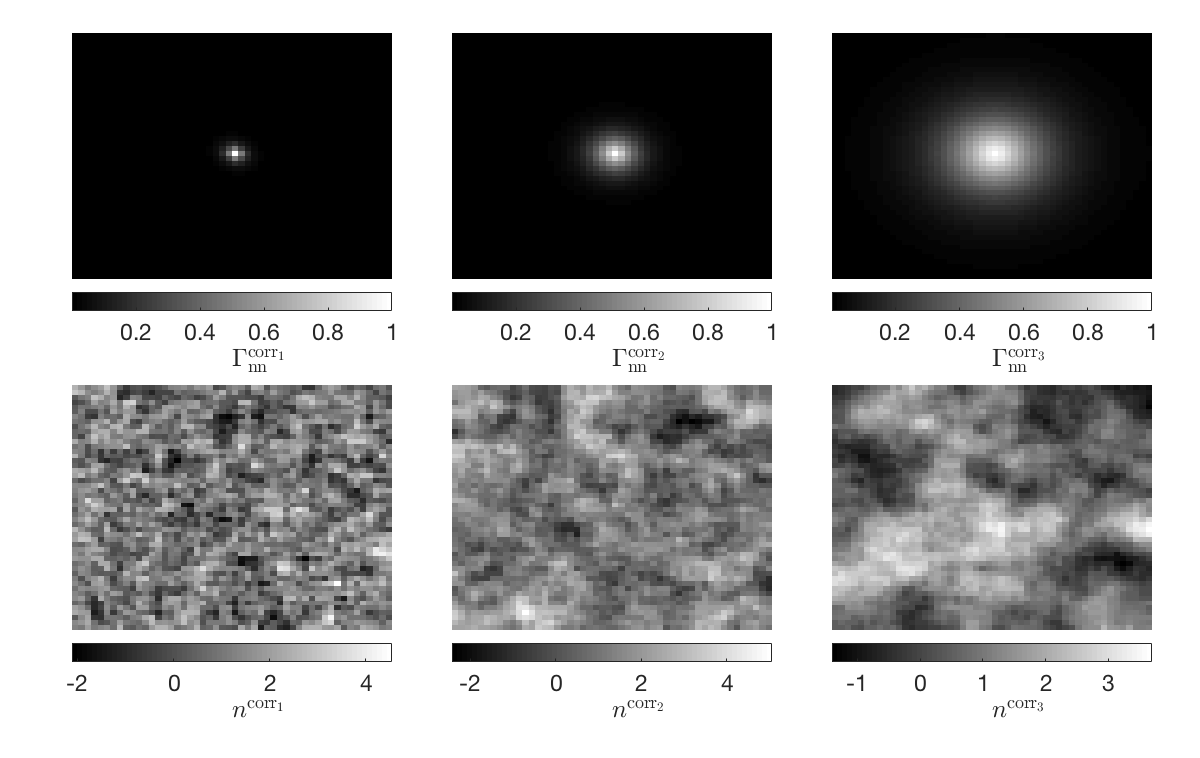} 
\caption{Top row.
Spatial correlation models for the multiplicative noise with different spatial decay rates.
Bottom row.
Draws from the respective models.}\label{fig:CorrNoise}
\end{figure}

\begin{figure}[h!]
\includegraphics[width=16cm]{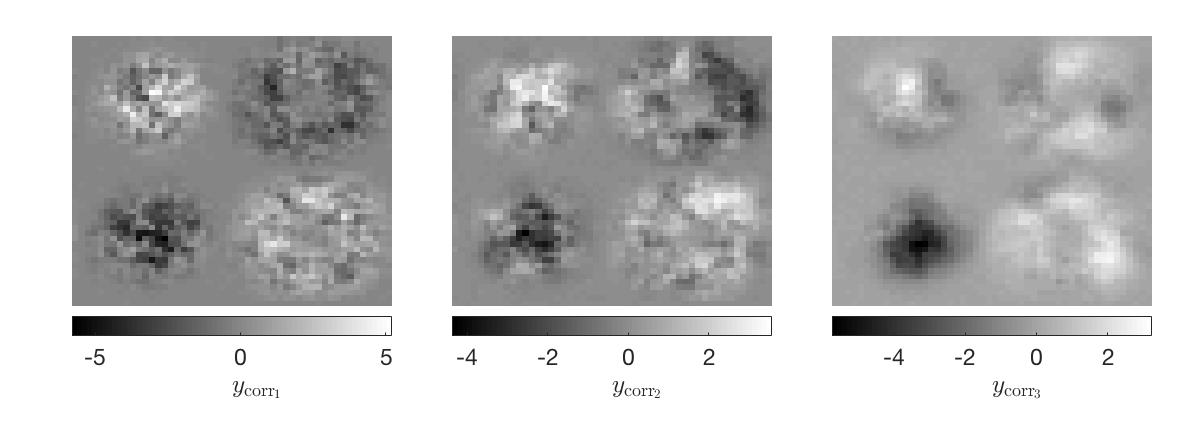} 
\caption{The observations with the three different spatially correlated multiplicative noise models
shown in Fig.~\ref{fig:CorrNoise}.}\label{fig:CorrNoiseData}
\end{figure}

In this section, we only consider normal multiplicative noise. 
We generate three distributions with different spatial decay rates. 
The traces of the multiplicative noise covariances are the same as in the cases of spatially uncorrelated
noise case.
Furthermore, the variance of the (spatially uncorrelated) additive noise is as in the previous case.
The respective correlation functions and draws from these distributions are shown in Fig.~\ref{fig:CorrNoise}.
The respective observations are shown in Fig.~\ref{fig:CorrNoiseData}.

\begin{figure}[h!]
\includegraphics[width=16cm]{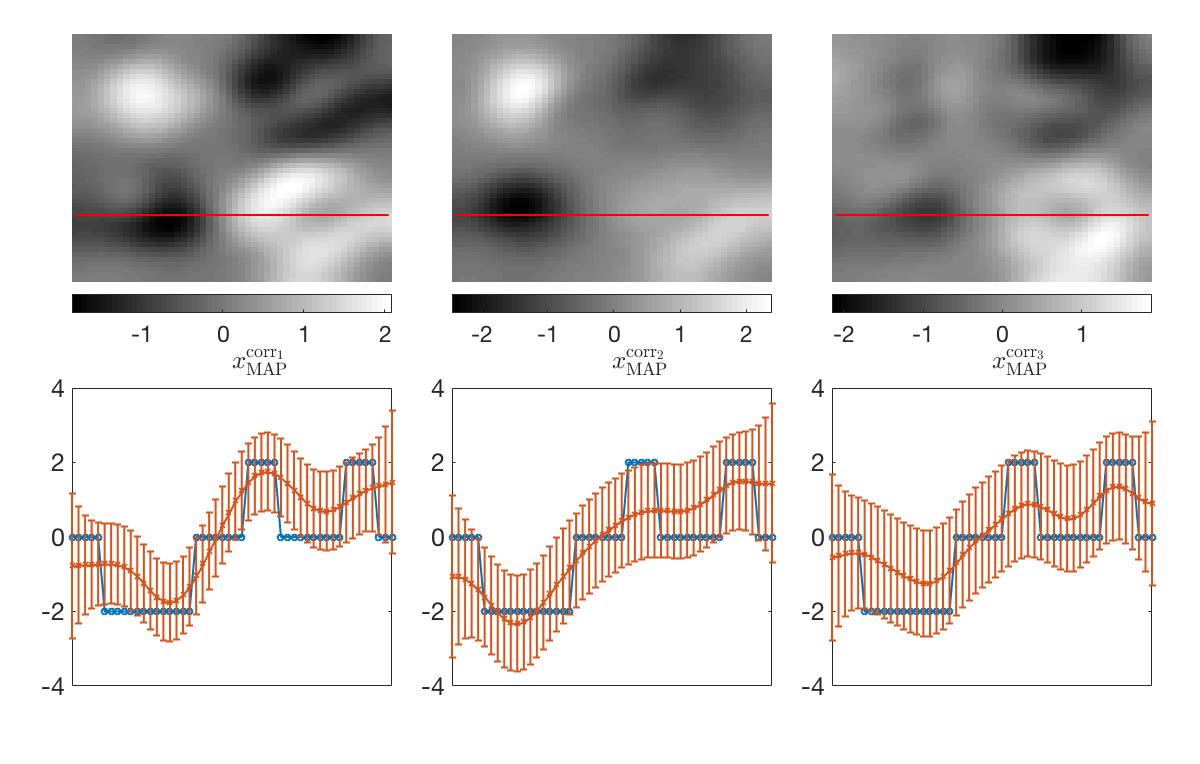} 
\caption{Top row.
The MAP estimates attained by using the BAE approach with different 
spatially correlated multiplicative noise models shown in Fig.~\ref{fig:CorrNoise}.
Left: Gamma, Centre: normal and Right: uniform noise models. 
Bottom row.
Cross sections of the actual target and reconstructions with  approximate MAP $\pm 3\sqrt{\bs{\Gamma}_{x\vert y}(k,k)}$
intervals along the lines in the top row reconstructions.}
\label{fig:CorrNoiseRes}
\end{figure}

The approximate MAP estimates and the posterior $\pm3$ STD intervals are shown in Fig.~\ref{fig:CorrNoiseRes}.
The estimates are, again, feasible with respect to the posterior error intervals.
The error estimates are larger than in the case of spatially uncorrelated multiplicative noise
which was expected.
As was also expected, the error estimates increase with decreasing decay rate of spatial correlation of the
multiplicative noise.

\section{Conclusion}

In this paper, we proposed an  approach to approximate (linear) inverse problems 
corrupted by both additive and multiplicative noise with an additive noise model.
The approximate additive noise model is constructed by (approximate) marginalization over the 
discrepancy between the model predictions of the original and the approximate model,
which is referred to as the Bayesian approximation error (BAE) approach. 
The resulting additive noise term is then approximated with a normal.
The covariance of this term is nontrivial and depends on the prior covariance.
The computation of the approximate MAP estimate does not suffer from convexity-related problems other
than those  arising from the forward map.

The approach does not need the multiplicative noise to be uncorrelated.
The results in this paper are, however, based on mutual independence of the primary unknown 
and the multiplicative noise.
As such, the mutual independence of the primary unknown and the multiplicative noise is, however, 
not essential for the proposed approach.
In such a case, the computation of the related joint covariance of the modified additive noise and the 
primary unknown involves rather tedious mappings of general fourth order statistics (kurtosis).

We considered numerical examples with different the multiplicative noise distributions related to 
an image processing deconvolution problem with both additive and multiplicative noise.
The results show that the approximation is feasible in the sense that the posterior error estimates support the 
actual target image.
Furthermore, the results were feasible also when the multiplicative noise was spatially highly correlated.

\dontshow{However, the approach outlined could be used for many other multiplicative noise models, such as Rayleigh distributed multiplicative noise. Moreover, for any distribution placed on the multiplicative noise the resulting functional to minimised in order to calculate the MAP estimate is quadratic, avoiding all of the convexity issues which are apparent with other methods equiped to deal with multiplicative noise.
Two distinct ideas have yet to be addressed. Firstly in some cases the product density referred to in Section \ref{sec: BAE} can be calculated in closed form, thus it may not be too difficult to drop the Gaussian approximation on the updated additive error term in favour of the retaining the true density. Secondly, the methods outlined in this paper do not require the presumption of a Gaussian prior to be placed on the parameter of interest, hence it would be of interest to apply these methods when using a TV-type prior density.
}

\bibliographystyle{siam}
\bibliography{MultNoise,jkadd}

\end{document}